\documentclass[a4paper]{amsart}
\usepackage{amscd}
\usepackage{amsmath}
\usepackage{amsfonts}
\usepackage{amssymb}
\usepackage{pstricks}
\usepackage{pst-node}

\numberwithin{equation}{section} \swapnumbers

 \DeclareMathOperator{\ke}{Ker}

 \DeclareMathOperator{\Spec}{Spec}

\def\+{{\dagger}}

\def\bd{{\begin{displaymath}}}
\def\ed{{\end{displaymath}}}

\def\arxiv{arXiv: \ignorespaces}

\def\AA{{\mathbb A}}
\def\CC{{\mathbb C}}
\def\DD{{\mathbb D}}
\def\FF{{\mathbb F}}
\def\GG{{\mathbb G}}

\def\PP{{\mathbb P}}
\def\QQ{{\mathbb Q}}
\def\RR{{\mathbb R}}
\def\ZZ{{\mathbb Z}}

\def\C{{\mathcal C}}
\def\D{{\mathcal D}}
\def\E{{\mathcal E}}
\def\F{{\mathcal F}}
\def\G{{\mathcal G}}
\def\cH{{\mathcal H}}
\def\I{{\mathcal I}}

\def\M{{\mathcal M}}

\def\cO{{\mathcal O}}
\def\cP{{\mathcal P}}

\def\U{{\mathcal U}}
\def\V{{\mathcal V}}

\def\X{{\mathcal X}}
\def\Y{{\mathcal Y}}

\newcommand{\fG}{\mathfrak{G}}

\newcommand{\fP}{\mathfrak{P}}

\newcommand{\fU}{\mathfrak{U}}

\newcommand{\fX}{\mathfrak{X}}

\newcommand{\fb}{\mathfrak{b}}

\newcommand{\fg}{\mathfrak{g}}
\newcommand{\fh}{\mathfrak{h}}

\newcommand{\fl}{\mathfrak{l}}
\newcommand{\fm}{\mathfrak{m}}

\newcommand{\fp}{\mathfrak{p}}

\newcommand{\fs}{\mathfrak{s}}
\newcommand{\ft}{\mathfrak{t}}

\newcommand{\fz}{\mathfrak{z}}

\begin{document}
\title[ ] {Arithmetic $\D$-modules and Representations }

\author[Lai]{King Fai Lai}
\address{Castle Peak Road, Tuen Mun\\
   Hong Kong}
\email{rexinterra@gmail.com}


\thanks{\textbf{Acknowledgement:} These are notes based on talks I
gave at the the Mathematics Institute of NCTS at Tsing Hua, Taiwan (on 9 May 2007),
 the Institute of Mathematics at the University Rennes (IRMAR) (on 29 June 2007)
 and the Mathematics Institute of Karlsruhe University (on 4 July 2007).
 Their generous hospitalities are gratefully
 acknowledged.}

\begin{abstract} 
We propose in this paper an approach to Breuil's conjecture on a
Langlands correspondence between $p$-adic Galois representations
 and representations of $p$-adic Lie groups in $p$-adic topological vector spaces.
We suggest that
 Berthelot's theory of arithmetic $D$-modules should give
a $p$-adic analogue of Kashiwara's theory of $D$-modules for real
Lie groups i.e. it should give a realization of the $p$-adic
 representations of a $p$-adic Lie group as spaces of overconvergent solutions of
 arithmetic $D$-modules which will come equipped with an action of the
 Galois group. We shall discuss the case of Siegel modular varieties
as a possible testing ground for the proposal.
\end{abstract}
\date{}

\maketitle

\section{Introduction}
\label{sec:intro}

Breuil conjectured that
   there is
 a Langlands correspondence between $p$-adic Galois representations
 and representations of $p$-adic Lie groups in $p$-adic topological vector spaces.
 Breuil-Schneider (\cite{BrS 03})
  pointed out that at the moment it is difficult to construct $p$-adic
 representations of a $p$-adic Lie group.

In this note I propose an approach to this problem, namely,  use Berthelot's
theory of arithmetic $D$-modules to give a $p$-adic analogue of
Kashiwara's theory of $D$-modules for real Lie groups i.e. we want
to realize $p$-adic
 representations of a $p$-adic Lie group as solution spaces of
 arithmetic $D$-modules which will come equipped with an action of the
 Galois group. This is described in the third section.
In the last section we consider the case of Siegel modular
varieties as a possible testing ground for the proposal.

As this note is based on talks to audiences with diverse backgrounds,
I mention facts which may be known to one group but not to another.
In any case this is a dream so no apologies to experts.

At this point it remains for me to thank Christoph Breuil, Mat Emerton,
Paul Gerardin, Michel Gros, Dick Gross, C-G. Schmidt, Jing Yu for conversations on this
note,
Pierre Berthelot  for
interest in and support for this fantasy, Raja Varadarajan for his help with the hypergeometric
equation and  Hotta Ryoshi
 for sending his Indian notes to an unknown foot soldier.

\section{In the beginning}

\subsection{Over $\RR$}\label{sec:sl2}
\smallskip

 Let me begin with
 the discrete series of $G = SL_2(\RR)$ (Bargmann \cite{Bar 47}).
The Lie algebra $\fg = {\fs \fl}_{2} (\RR)$ of $SL_{2}(\RR)$
consists of $2 \times 2$ matrices of trace 0 and is generated by
\[ X = \left( \begin{array}{ll}
0 & 1 \\
0 & 0
\end{array}\right),  \hspace{5mm}
Y = \left( \begin{array}{ll}
0 & 0 \\
1 & 0
\end{array}\right) \hspace{5mm} \mbox{and} \hspace{5mm}
H = \left( \begin{array}{ll}
1 & 0 \\
0 & -1
\end{array}\right) \, .\]
The center of the universal enveloping algebra $U(\fg)$ is generated
by the Casmir operator :
\[\Omega = (H+1)^2 + 4 YX
\]

We identify elements of the universal enveloping algebra $U(\fg)$
with differential operators on $G$. If we use polar coordinates via the
Cartan decomposition of $G$:
\[ \phi : (\theta_1, t,   \theta_2) \mapsto u_{\theta_1} a_t  u_{\theta_2} \]
where $u_{\theta} = \left( \begin{array}{ll}
\cos \theta & \sin \theta \\
-\sin \theta & \cos \theta
\end{array}\right)$ and $a_t =
\left( \begin{array}{ll}
e^t & 0\\
0 & e^{-t}
\end{array}\right)$
then $d\phi$ is bijective and a differential operator $D$ on $G$
goes to the differential operator $\tilde{D}$. For $\Omega$ this
gives
\[
\tilde{\Omega} = \frac{1}{\sinh^2 2t}\left(
\frac{\partial^2}{\partial \theta_1^2} + \frac{\partial^2}{\partial
\theta_2^2} \right) -2\frac{\cosh 2t}{\sinh ^2 2t}
\frac{\partial^2}{\partial \theta_1 \partial \theta_2} +
\frac{\partial^2}{\partial t^2} + 2\frac{\cosh 2t}{\sinh
2t}\frac{\partial}{\partial t} +1
\]

Let $k \geq 1$ and $\pi_k$ be the discrete series with lowest weight
$k+1$ and infinitesimal character $\chi_k$ i.e. $\chi_k (\Omega) =
k^2$. We can realize this representation in $L^2(G)$.
 Let us take a square integrable $f$ satisfying
\begin{equation*} f(u_{\theta_1}xu_{\theta_2}) = e^{i(m \theta_1 + n \theta_2)}f(x).
\end{equation*}
Since
\begin{equation*}
dg = \frac{1}{2} \sinh 2t d\theta_1 dt d\theta_2
\end{equation*}
we put $ \tilde{f} (t) = (\sinh 2t )^{\frac{1}{2}} f(a_t)$. Then in
polar coordinates the Casmir operator $\tilde{\Omega}$ goes to the
operator $\Delta_{mn}= d^2 /dt^2 - q_{mn}$ with
\begin{equation*} q_{mn} = \frac{m^2 +n^2 - 2mn \cosh 2t -1}{\sinh ^2 2t}.
\end{equation*}
Then function $\tilde{f}$ then satisfies the differential equation
\begin{equation*} \Delta_{mn} \tilde{f} = k^2 \tilde{f}.
\end{equation*}
If we go to a new variable $u= \cosh 2t$ then the differential
equation becomes
\begin{equation*}
\frac{d^2 g}{du^2} + \frac{P(u)}{u^2 -1} \frac{dg}{du}+
\frac{Q(u)}{(u^2 -1)^2}g = 0
\end{equation*}
where $P, Q$ are polynomials in $u$ and $ \deg \;P \leq 1$, $ \deg
\;Q \leq 2$. This is a hypergeometric equation with singularities at
$u=1, -1, \infty$ - such equations are indeed very familiar in
$p$-adic analysis which has been studied by Dwork and his school
(\cite{Dwo 73}, \cite{Put 86}, \cite{Ked 05}).
Notice there is one equation for each pair $(m,n)$ i.e. for each $K$
type - this will be important when we put coherence conditions on admissibility to insure
some holonomy. In addition for something closer to modular forms
I would like to draw attention to the Halphen Fricke differential operators and
its application to $p$-adic modular forms in \cite{Kat 76} (in particular \S 2.4.5).
For more differential operators on automorphic forms
see \cite{Shi 81} and its applications to $p$-adic Siegel
modular forms see \cite{Pan 04}.

We have the usual action of  $G$  on the upper half plane $\fh$ by :
\[ g \cdot z = \frac{az+b}{cz+d} \, \]
with $g = \left(
\begin{array}{ll}
a & b \\
c & d
\end{array}\right) \in G$ and $z \in \fh$.
The discrete series representation $\pi_k$ can also be realized in the space of
holomorphic functions $f: \fh \to \CC$ on the upper plane $\fh$
which are square integrable with respect to the measure $y^{k-2} dx
dy$. The action of $SL_2(\RR)$ is given by
\begin{equation*} \left( \pi_k
\left(
\begin{matrix}
a & b \\
c & d
\end{matrix}
\right) f \right) (z) = (-cz +a )^{-k} f (\frac{dz-b}{-cz+a})
\end{equation*}
for $z \in \fh$. If we put
\[ e_n (z) = (z-i)^{n-k-1} (z+i)^{1-n} \; ,\]
then the space of this representation is the completed Hilbert sum
of $e_n$ for $n \geq k+1$. This is often referred to as the holomorphic discrete series.

\subsection{The $p$ adic case}
Next we consider the $p$-adic representations of Berger-Breuil. Let
$E$ be an extension of $\QQ_p$. For a positive real number $r$, say
a function $f: \ZZ_p \to E$ is in $\C^r$ if $n^r|a_n(f)| \to 0$ in
$\RR$ where $a_n(f) = \sum^n_{i=0} (-1)^i \left(
\begin{array}{l}
n \\
i
\end{array}\right) f(n-i)$.
Note that $f(z) = \sum_{n=0}^{\infty} a_n(f)\left(
\begin{array}{l}
z \\
n
\end{array}\right)$ and that $f$ is in $\C^n$ if and only if $f$ has
{\it continuous} $p$-adic derivative of order up to $n$.

Normalize the valuation by $val(p) =1$. Write the ring of integers
of $E$ as $\cO_E$. Take $\alpha \neq \beta$ in $\cO_E$. Let
$B(\alpha)$ be the space of functions $f: \QQ_p \to E$ such that the
restrictions $f|\ZZ_p$ is in $\C^{val(\alpha)}$ and $(\alpha p
\beta^{-1})^{val(z)}z^{k-2}f(1/z) | \ZZ_p$ can be extended to a
function in $\C^{val(\alpha)}$. Then $GL_2(\QQ_p)$ acts in the
following manner:
\begin{equation*}
\left(
\begin{matrix}
a & b \\
c & d
\end{matrix}
\right) f  (z) = \alpha^{val(ad-bc)}(\alpha p
\beta^{-1})^{val(-cz+a)} (-cz +a )^{k-2} f (\frac{dz-b}{-cz+a}).
\end{equation*}
Let $L(\alpha)$ be the closed subspace of $B(\alpha)$ spanned by
$z^j$, $(\alpha p \beta^{-1})^{val(z-a)} (z  - a )^{k-2-j}$ for $a
\in \QQ_p$, $j \in \ZZ$, $0 \leq j < val(\alpha).$ Put $\Pi =
B(\alpha)/ L(\alpha).$ Then it is a deep theorem of Berger-Breuil
that $\Pi \neq 0$ and that $\Pi$ is an admissible representation od
$GL_2(\QQ_p)$.

Looking at the action defined by Berger-Breuil and the action in the
holomorphic discrete series of $SL_2(\RR)$ one could not help asking
if one could write down the differential equation in this case
 or if it is possible to construct their model using
arithmetic theory of $\D$-modules. In fact
there is a differential equation lurking behind the model of
Berger-Breuil; in  \cite{BBr 06} \S 5.2
using the Fontaine functor (\cite{Fon 90} A 3.4) the dual of the representation $\Pi$ above
is shown to be isomorphic to the
space of bounded sequences in the inverse limit of $(\phi, \Gamma)$-modules
constructed from the starting data $\alpha, \beta$ and the inverse of Berger functor
would take a $(\phi, \Gamma)$-module to the solution space of a differential equation
(\cite{Beg 04} Theorem A and C).
The study of the differential equation on the boundary of the unit dics
reminds me of the asymptotics at infinity of invariant eigendistributions of Harish-Chandra.

\subsection{Philosophy}
Before I leave $SL_2$ I turn to a bit of general background.
The theory of Langlands correspondence says there is a bijection
between the representations $\pi$ of a reductive group $G$ over a
field $F$
 into the automorphism group
  of a vector space over a field $E$ and the
 representations $\rho$ of a Galois group of the field $F$ into the
 Langlands dual of  $G$
which "almost" matches the Frobenius eigenvalues of $\rho$ with the
spherical Hecke
 eigenvalues of $\pi$.
This has been established
 in the cases of $G= GL_n$, $E= \CC$  and $F$ is a
 a local field field,
 a one variable function field over a finite field or over the complex numbers
 by Langlands-Deligne-Kutzko-Harris-Taylor-Henniart, \;\;
 Laumon-Rapoport-Stuhler, \; \;
 Drinfeld-Lafforgue, \; \; Frenkel-Gaitsgory-Vilonen respectively
 (\cite{Car 00}, \cite{LRS 93}, \cite{Lau 02},\cite{Fre 04})
 and almost for $G=GL_2$, $E=\CC$, $F$ a totally real filed by
 Langlands-Deligne-Carayol-Saito.
(Readers will object to the lack of precision of the description
here -
in case $F$ is a number field one should say $\pi$ is an automorphic
representation and $\rho$ is a compatible system of $\lambda$ adic representations; or
one should mention the Weil-Deligne group, the automorphic
Langlands group (\cite{Art 02}), the motivic Galois group; \; \; one
should be more precise about the representations occurring and
mention automorphic forms, the L-indistinguishability of
Labesse-Langlands, endoscopy of Kottwitz-Shelstad; \; \; one should
point out the Deligne-Laumon geometric interpretation with
 Galois group as fundamental group and using Grothendieck's function-sheaf dictionary
 on the bijection between $GL_n(F) \backslash  GL_n(\AA) / GL_n(\cO)$
 isomorphism classes of rank $n$
 vector bundles.
But this would be besides the point here.)

Reading the  proofs of the known correspondences one could not help
having a feeling of staring at a magical inventory  in a gigantic
Amazon warehouse. The purpose of this note is to propose a "natural"
setting for the Langlands correspondence, namely, the Galois group
and the reductive group occurs as mutual centralizers in the
monodromy group of a system of differential equations. This might
not sound so far fetched from the points of view of Grothendieck's
theory of monodromy and Harish-Chandra's constructions of
representations using differential equations.

\subsection{Harish-Chandra}
Let $G$ be a semisimple Lie group with finite center and $\Gamma$ be
a discrete subgroup of $G$ such that $\Gamma \backslash G$ has
finite volume. Then $G$ acts on the Hilbert space of square
integrable functions $L^2(\Gamma \backslash G)$ by translation $f(x)
\mapsto f(g^{-1}x)$ giving the so called the regular representation
of $G$.
 A basic theme in the theory of automorphic forms
is to study the representations of $G$ that occur in the  spectral
decomposition of this regular representation. For this we need to
know the representations of $G$. This is done by Harish-Chandra
using the theory of differential equations. Differential equations
comes in right from the beginning in the paper of Bargmann which was
extended by Harish-Chandra to all bounded symmetric domains, this
line leads eventually to Schmidt's solution of Langlands conjecture
constructing representations on $L^2$-cohomology of line bundles.

One more remark to make is this: following the tradition of
Local-Global principle in Classfield theory, we replace $\Gamma
\backslash G$ by $G(\QQ)\G(\AA)$ for a reductive  group $G$ over
$\QQ$ and study representations of the form $\otimes \pi_p$ with
$\pi_p$ a representation of $G(\QQ_p)$ on complex vector spaces
(here: for $p=\infty$ we take $\QQ_p$ to be $\RR$). Again
Harish-Chandra has constructed these representations for finite $p$.

\subsection{D-modules}
Nowadays beginning with Beilinson-Bernstein we discuss the
representation theory of real Lie groups in the frame work of
D-modules.

For a complex analytic (or algebraic) manifold $X$ let
$\D_X$ be the sheaf of analytic differential operators on $X$
(\cite{Bjo 93} \S 1.2.2; or \cite{Bor 87} II 2.12, VI 1.1, 1.2).
On an open subset $U$ of $X$,
to an differential operator $\Delta$ in $\D_X(U)$ we associate
the $\D_X(U)$-module $M + \D_X(U) / \Delta \D_X(U)$.
Then we can think of $Hom_{\D_X(U)}
(M, \O_X(U))$ as the space of solutions of the differential equations
$\Delta = 0$ on $U$. This will be our object of interest.

Partly to
supplement the claim  that "representations of a
semisimple Lie group are realized as solution spaces of systems of
differential equations"  I recall some results of
Kashiwara - Schmid company.

Let $G$ be a semisimple Lie group with
finite center with Lie algebra $\fg$. Let $K$ be the maximal compact
subgroup of $G$. Let $Mod (\fg_{\CC}, K)$ denote the category of
algebraic $(\fg_{\CC}, K)$-modules; this category contains the
Harish-Chandra modules which correspond to admissible
representations of $G$. Write $D^b (Mod (\fg_{\CC}, K))$ for the
bounded derived category of $Mod (\fg_{\CC}, K)$. Let $Mod_G
(\D_{G/K})$ denote the category of quasi coherent $G$ equivariant
$\D_{G/K}$-modules. Write $D^b_G (\D_{G/K})$ for the bounded derived
category of $Mod_G (\D_{G/K})$. The Kashiwara showed that the
categories $D^b (Mod (\fg_{\CC}, K))$ and $D^b_G (\D_{G/K})$ are
equivalent.

Again let $G$ be a semisimple Lie group with finite center with Lie
algebra $\fg$. Write $\U (\fg)$ for the universal enveloping algebra
of $\fg$. Fix a Borel subgroup $B$ of $G$. Put $\U_{G/B}({\fg}) =
\cO_{\G/B} \otimes \U (\fg)$ with multiplication such that $\fg$
acts on $\U (\fg)$ by left multiplication and on $\cO_{\G/B}$ by
differentiation. Let $\tilde{\fg}$ be the kernel of the morphism
from $\U_{G/B}({\fg})$ to the sheaf of vector fields on $G/B$. Fix a
maximal torus $T$ in $B$. Write the Lie algebra of $T$ as $\ft$. Let
$\lambda$ be a $\CC$-algebra homomorphism from the center $\fz$ of
$\U (\fg)$  to $\CC$ which we also identify with an element of the
dual $\ft^*$ and extend to an element of $\fb^*$. We get a one
dimensional representation $\CC_{\lambda} = \CC \cdot 1_{\lambda}$
of $\fb^*$ by $A \cdot 1_{\lambda} = \lambda(A) \; 1_{\lambda}$ for
$A \in \fb$. Set
$$
\underline{A}_{G/B} (\lambda) = \U_{G/B}({\fg}) / \sum_{A \in
\tilde{\fg}} \U_{G/B}({\fg}) (A - \lambda(A)).
$$
Let $\rho$ be half sum of the positive roots of $(\fg_{\CC},
\ft_{\CC})$ with respect to $\fb$. The ring $\D_{G/B, \lambda} =
\underline{A}_{G/B} (\lambda + \rho)$ is Kashiwara's twisted ring of
differential operators on $G/B$ and he proved in this case a
Riemann-Hilbert correspondence which says that $\M \mapsto
R{\cH}om_{\D_{G/B, \lambda}} (\cO_{G/B} \otimes \CC_{\lambda}, \M)$
yields an equivalence from the bounded derived category of regular
holonomic $ \D_{G/B, \lambda}$-modules to the bounded derived
category of constructible sheafs of complex vector spaces on $G/B$
twisted by $-\lambda$.

Before I forget I would like to mention the work of Hotta (\cite{Hot 87}) on
character sheaves ( and the related work in the $p$ adic case of Gros \cite{Gro 01},
Aubert-Cuningham \cite{AuC 04} ) and on Morihihiko Saito ( and the related comments of
Berthelot in the case of arithmetic modules).


\section{Homogeneous Formal Schemes}

\subsection{Equivariant Arithmetic $\D$-modules}\label{sec:dmod}

Fix a prime $p$. Let $\V$ be a complete discrete valuation ring of
characteristic zero with maximal ideal $\fm$ and residue field $\FF$
of characteristic $p$. Write $F$ for the field of fraction of $\V$.

Let $\fX$ be a locally noetherian proper smooth formal $\V$ scheme.
we always assume that ${\fm}{\cO}_{\fX}$ is an ideal of definition.

Let $\fG$ be an affine formal group $\V$ scheme acting homogeneously
on $\fX$ over $\V$. Let us write the action map as $m: \fG \times
\fX \to \fX : g, x \mapsto gx$ and the projection as $p: \fG \times
\fX \to \fX$. Define maps $d_j : \fG \times  \fG \times \fX \to \fG
\times \fX$ by $d_1 (g_1, g_2, x) =(g_1, g_2 x)$, $d_2(g_1, g_2, x)
=(g_1 g_2, x)$, $d_3(g_1, g_2, x) =(g_2, x)$. Then our assumption
says $m d_1 = m d_2$, $p d_2 = p d_3$, $m d_3 = p d_1$. To begin
with this action would allow us to consider the elements of the
universal enveloping algebra of $\fG$ as differential operators on
$\fX$.

Let ${\D}_{\fX}^{\+}$ the sheaf of differential operators of
infinite order of finite level on $\fX$ introduced by Berthelot
(\cite{Ber 961} (2.4.1)). By a $\fG$-equivariant ${\D}_{\fX}^{\+}$
module we mean a pair consisting a ${\D}_{\fX}^{\+}$ module $\M$
together with a ${\D}_{\fG \times \fX}^{\+}$ linear isomorphism $b:
m^* \M \to p^* \M$ such that we have the following commutative
diagram
\begin{equation*}
\begin{CD}
 d_2^* m^* \M @>d_2^* b>> d_2^* p^* \M\\
@VVV @VVV  \\
d_1^* m^* \M @>a>> d_3^* p^* \M
\end{CD}
\end{equation*}
where $a$ is the composite $d_1^* m^* \M \xrightarrow{d_1^* b}
d_1^* p^* \M \to d_3^* m^* \M \xrightarrow{d_3^* b}
 d_3^* p^* \M$.

Let $Mod_{\fG} ({\D}_{\fX}^{\+})$ denote the category of {\it
admissible} $\fG$-equivariant ${\D}_{\fX}^{\+}$ modules with $\fG$
equivariant morphisms between them. We have left open the
requirements that should be put into the word {\it admissible}. For
example we would want $Mod_{\fG} ({\D}_{\fX}^{\+})$ to be an abelian
category so that we can define the associated bounded derived
category $D_{\fG}^b ({\D}_{\fX}^{\+})$  of $\fG$-equivariant
${\D}_{\fX}^{\+}$ modules. Another point is the condition of
$K$-coadmissible for compact subgroup $K$ in the definition of an
admissible $p$-adic representation should be translated into some
"coherence" condition on the ${\D}_{\fX}^{\+}$-module.

\subsection{Equivariant categories}

Let $\X$ be a (quasi-separated) rigid analytic space whose formal
model is a proper smooth formal $\V$ scheme (\cite{BLR 93}). Fix a
complete valued field $E$ containing ${\QQ}_p$. Assume that
$\cO_{\X}$ is an  $E$ algebra. The category of sheaves of $\cO_{\X}$
modules which are locally convex vector spaces of compact type over
$E$ is exact; let $D^b(\X)$ be its bounded triangulated category.
(Recall that a p-adic locally convex topological space is said to be
of compact type if it is a direct limit of Banach spaces with
compact transition maps.) Let $F$ be the field of fractions of $\V$.
Let $G$ be the $F$ rational points of a reductive group over $F$.
Assume that $G$ acts on $\X$. Then a $G$ equivariant $E$ sheaf is an
$\cO_{\X}$ module $M$ of locally convex vector spaces of compact
type over $E$ together with an $\cO_{\X}$ linear isomorphism $b: m^*
M \to p^* M$ satisfying the usual cocycle condition given by a
commutative diagram as above. The category $D^b(\X)$
 contains the bounded triangulated category of $G$ equivariant $E$ sheaves.
Let $D^b_G(\X)$ be the $G$-equivariant  bounded triangulated
category following Bernstein-Lunts (according to Gottker-Schnetmann)
and let $\F : D^b_G(\X) \to D^b (\X)$ be the forgetful functor such
that (1) under $\F$ the image of the heart of $D^b_G(\X)$ contains
the $G$ equivariant sheaves and (2) $D^b_G(\X)$ is equal to $D^b( G
\backslash \X)$ if the quotient $G \backslash  \X$ exists as a rigid
space.

\subsection{Overconvergent Solutions}\label{sec:solu}

\subsubsection{}
As I have taken the liberty to write out the formulas for
$SL_2(\RR)$ I might as well do the same for "overconvergence" for
easy reading. Fix a prime $p$. Let $\V$ be a complete discrete
valuation ring of characteristic zero with maximal ideal $\fm$ and
residue field $\FF$ of characteristic $p$. Write $F$ for the field
of fraction of $\V$. Let $\fP$ be a locally noetherian proper
formal $\V$ scheme. Write $P_{\FF}$ for its closed fibre, $\cP$ for
the rigid analytic space associated to its generic fibre.

Let $sp: \cP \to \fP$ be the specialization
 map according to  Berthelot (\cite{Ber 96} (0.2.2)) -
locally it is the reduction modulo $\fm$. Let $j: U \to X$ be an
open immersion of $\FF$-varieties and $X \to P_{\FF}$ be a closed
immersion. For example $X = P_{\FF}$. We write this data as a frame
$U \xrightarrow{j} X \to \fP$. We write $]U[$ for the inverse image
$sp^{-1} (U)$ which we shall call the tube of $U$. We say a
neighbourhood $V$ of $]U[$ in $]X[$ is strict if $]X[ = V \cup ]X -
U[$ is an admissible covering in the rigid topology (\cite{Ber 96}
(1.2.1)).

For a sheaf $\E$ on an admissible open subset $V$ of $]X[$ we define
the sheaf of germs of sections of $\E$ overconvergent along $X - U$
or simply the overconvergent sheaf
\[ j^{\+} \E = {j_{V}}_{*} \lim_{\substack{\longrightarrow}} {j_{VV'}}_{*} {j_{V V'}}^{*} \E
\]
where $V'$ runs through all the strict neighbouhoods of $]U[$ in
$]X[$, $j_{VV'} : V \cap V' \to V$ and $j_{V} : V \to ]X[$
(\cite{Ber 96} (2.1.1)).

To look at the affine picture let us write $t$ for $(t_1, \cdots,
t_n)$. Then the weak completion of the polynomial ring is
\[ \V[t]^{\+} = \{\sum a_{\alpha} t^{\alpha} \in \V[[t]] :
\text{for some} \; r > 1 , |a_{\alpha}|r^{|\alpha|}  \to 0 \} \; ,
\]
and the $\fm$-completion of $\V[t]^{\+}$ is
\[\V<t> = \{\sum a_{\alpha} t^{\alpha} \in \V[[t]] :
\lim a_{\alpha} = 0 \} \;.
\]
Here we use multi-index notation for $t^{\alpha}$. Now take $A=
\V[t] / (f_1, \cdots, f_m)$, $U = \Spec A$. Let $X$ be the closure
of $U$ in the projective space ${\PP}^n_{\V}$, $\fX$ be the
completion along the closed fibre. We have a frame $(U_{\FF} \subset
X_{\FF} \subset \fX)$. Let us write $B^n(0,r^+)$ for the closed disc
of radius $r$ and $\U$ for the rigid analytic space associated to $U
\otimes_{\V} F$. The admissible open sets
\[ V_r = B^n(0,r^+) \cap \U
\]
form a cofinal family of affinoid strict neighbourhoods of the tube
$]U-{\FF}[$ in $]X_{\FF}[$ for $ r \geq 1$, $r \to +1$. Put $A_r =
\Gamma(V_r , \cO_{\V_r}).$ Write
$$F<\frac{t}{r}> = \{ \sum b_{\alpha} t^{\alpha} \in F[[t]] :
\lim b_{\alpha}r^{|\alpha|} \to 0 \} \;.$$
 Then
 $$A_r = F<\frac{t}{r}> / (f_1, \cdots, f_m)$$
 and
$$(A \otimes_{\V} F )^{\+} = \lim_{\substack{ r \to +1}} A_r =
\V[t]^{\+}/(f_1, \cdots, f_m) \otimes_{\V} F    \;.$$ Moreover the
category of coherent $j^{\+} \cO_{\U} $-modules is the same as the
category of coherent  $(A \otimes_{\V} F )^{\+}$ - modules . I hope
that this helps.

\subsubsection{}
We continue in the notation of section \ref{sec:dmod}. Given a
formal $\V$ scheme $\fX$ we denote the rigid analytic space
associated to the generic fibre of $\fX$ by $\fX_F$ or simply $\X$.
To an $\cO_{\fX}$ module $M$ we can associate an $\cO_{\X}$ module
which we denote by $M_F$. Let $sp: \X \to \fX$ be the specialization
map of Berthelot. Let $j^{\+} \cO_{\X}$ be the sheaf of
overconvergent "functions" on $\X$.

Let $\M$ be in the bounded derived category $D_{\fG}^b
({\D}_{\fX}^{\+})$  of $\fG$-equivariant ${\D}_{\fX}^{\+}$ modules.
Then  applying the solution functor and then taking the associate
sheaf over the rigid fibre we get $R{\cH}om_{{\D}_{\fX}^{\+}} (\M,
sp_{*} j^{\+} \cO_{\X} )_F$ lying in $D_G(\X)$. Applying the global
section functor $\Gamma$ yields representations of $G$
on locally convex topological vector spaces over the $p$ adic fields
which we
simply denote by $Ext^*(\M, \cO_{\X})$. The aim is to find the
conditions for "{\it admissibility}" for $\M$ so that $Ext^*(\M,
\cO_{\X})$ are admissible representations in the sense of
Schneider-Teitelbaum. If this the case and as the Galois group
$\Gamma$ of the algebraic closure of $F$ over $F$ is acting
everywhere, we can try to decompose $Ext^*(\M, \cO_{\X})$ as an $G
\times \Gamma$ space and so obtain a natural setting for pairing
$p$-adic Galois representations with $p$-adic representations of the
$p$-adic reductive group $G$.

\subsection{Some definitions}
\subsubsection{}
It will be convenient to recall the definition of admissible
representation here.

Fix a complete $p$-adic field $E$. Given a {\it compact} $p$-adic
Lie group $K$ let ${\DD}(K)$ denote the strong dual of the space of
$E$-valued analytic functions on $K$. Then ${\DD}(K)$ is an inverse
limit of noetherian Banach algebras $A_n$ with {\it flat} transition
homomorphisms. Say a ${\DD}(K)$ module $W$ is $K$-coadmissible if it
is isomorphic to an inverse limit of {\it finitely generated} $A_n$
modules $W_n$ which are equipped with isomorphisms $A_n
\otimes_{A_{n+1} } M_{n+1} \to M_n$.

A representation of a $p$-adic group G on a locally convex space $V$
over a complete $p$-adic field $E$ is said to be locally analytic if
$V$ is barrelled and for all $v \in V$, the map $G \to V$ $:g
\mapsto gv$ is in the space of locally analytic maps from $G$ to
$V$. (For the definition of locally analytic maps from a $p$-adic
manifold to a $p$-adic locally convex space,
see \cite{ST1 02} section 2.) A locally analytic representation
of a $p$-adic group G on a locally convex space $V$ is said to be
{\it admissible} if $V$ is of compact type and the strong dual of
$V$ is $K$-coadmissible for all compact subgroups $K$ of $G$.

When the compact group $K$ is an union of affinoid open compact subgroups
of radii going 1, Emerton (\cite{Eme 04} Prop 5.2.3) identifies
the dual of the space of rigid analytic functions
on $K$ with a direct limit of $p$ adic completion of finite level
divided powers ( as in Berthelot \cite{Ber 02}) of the universal
enveloping algebra of the Lie algebra of $G$. This allows to show that
$\DD(K)$ is the inverse limit of rigid analytic distributions
(page 100 end of the proof of Cor 5.3.19) and eventually gives a definition of
admissibility equivalent to that of Schneider-Teitelbaum
but closer to the theory of arithmetic $D$ modules (\cite{Eme 04} definition6.1.1, Theorem 6.1.20).

\subsubsection{}
I would like to give the definition of the arithmetical differential operators.
I agree to it here is not in the logical order but fairy tales never is, besides
I feel that it is better to say it here after telling the theme for at this point
you will  see the familiar features. Before I begin I would also like
to recommend to the readers the wonderful papers of Grothendieck
(\cite{Gro 66}, \cite{Gro 74}, Chap IV).

Next we introduce divider powers of level $m$ following Berthelot
{\cite{Ber 02}). Write $\ZZ_{(p)}$ for the ring of integers localized at
the prime $(p)$. Let $A$ be a $\ZZ_{(p)}$-algebra.
For a given an ideal $J$ of $A$ a divided power structure on $ J$ is
a set of maps $\\gamma_n : J \to A \}$ satisfying certain axioms such that
$n! \gamma_n (x) = x^n$ for $ x \in J$ (see \cite{Gro 74} IV.1 for the full definitions).
For example take $A$ to be the quotient ring
$\ZZ / p^m \ZZ$ for an integer $m >1$, $J$ to be the ideal generated by
$p$, then $\gamma_n (x) = \frac{x^n}{n!}, x \in J$ defines a divided power structure
because $\nu_p (\frac{p^n}{n!}) \geq 1$.

A divided power structure of level $m$ is the data consisting of a $\ZZ_{(p)}$-algebra
$A$, an ideal $J$ equipped with divided power $\\gamma_n : J \to A \}$ and
an ideal $I$ of $A$ such that $I^{(p^m)} + pI \subset J$. Here we write
$I^{(h)}$ for the ideal generated by $x^h$ with $x \in I$. Often say
$I$ is the $m$-PD ideal. For example if
for all $x \in I$ we have $x^{p^m} \in pA$ then the ideal $J = pA \cap I$ with the canonical divided
power $\gamma_n (x) = \frac{x^n}{n!}$ would give a
divided power structure of level $m$.

Write a positive integer $k$
as $k = r +  p^m q$ with $0 \leq r < p^m$, and write $x^{[n]}$ for $\gamma_n (x)$,
introduce an operation
$$x^{\{ k\}} =x^r (x^{p^m})^{[q]}.
$$
The operation $x \mapsto x^{\{ k\}}$ enjoys analogous properties as divided powers.
Extend the divided power $x^{[n]}$ to $J_1 = J + (p)A$ (and use the same notation)
We shall write $I_k$ for the ideal of $A$ generated by
$x_1^{\{ n_1 \}} \cdots  x_t^{\{ n_t \}}$ for
$n_1 + cdots + n_t \geq k$ and write
$I^{\{n \}}$ for the ideal of $A$ generated by
$x_1^{\{ n_1 \}} \cdots  x_r^{\{ n_r \}}\; y_1^{[q_1]} \cdots  y_s^{[q_s]}$
with $x_i \in I$, $y_j \in (J +(p)A) \cap I_{k_j}$,and
$n_1 + \cdots + n_r + k_1 q_1 + \cdots + k_s q_s \geq n$.

Given a
$\ZZ_{(p)}$-algebra $A$
and any ideal $I$ in $A$ there exists an $A$-algebra
$P_{(m)} (I)$ equipped with  a divided power structure $(\bar{I}, \bar{J}, \bar{\gamma})$ of level $m$
satisfying
$IP_{(m)} \subset \bar{I}$
which are universal for the homomorphisms from $A$ to $A'$ carrying $I$ into
the $m$-PD ideal of $I'$ of $A'$.
We write
$$P_{(m)}^n (I) = P_{(m)} (I)/ {\bar{I}}^{\{n+1 \}}.
$$

Let $S$ be a scheme over $\ZZ/ p^{m+1} \ZZ$ and $X$ be a smooth
$S$ scheme. Let $\I$ be the ideal  in $\cO_{X \times X}$ of the
diagonal immersion $X \to X \times_{S} X$. The construction of
$P_{(m)}^n (I)$ leads to the construction of the sheaf
$\cP_{(m)}^n (\I)$ of principal parts of level $m$ of order $n$.
Let $\cO_{X}$ act on it via $f \mapsto f \otimes 1$ (for $f$ in
$\cO_{X}$). We define the sheaf of differential operators of level
$m$ of order $n$ as
$$\D^{(m)}_{X , n}
:= {\cH}om_{\cO_{X}} (\cP_{(m)}^n (\I),\cO_{X})
$$
and the sheaf of differential operators of level $m$
is $\D^{(m)}_{X } := \cup_{n \geq 0} \D^{(m)}_{X , n}.$

Now let  $\fX$ be  a smooth formal scheme  over
 a complete discrete valuation ring $\V$ with maximal ideal $\fm$.
Put $S_i = Spec (\V / \fm^{i+1})$ and $X_i = S_i \times_{\V} \fX$.
Define the sheaf of differential operators of level $m$ on the formal scheme $\fX$
by
$$\widehat{\D}^{(m)}_{\fX } = \lim_{\leftarrow \atop i} \D^{(m)}_{X_i }
$$
and finally the sheaf of arithmetical differential operators  on the formal scheme $\fX$
by
$$\D^{\dagger}_{\fX } := \cup_{m \geq 0} \widehat{\D}^{(m)}_{\fX }.$$
On an affine open subset $\fU$ of $\fX$ an element of
$\Gamma (\fU, \D^{\dagger}_{\fX })$ can always be  written as
$\sum a_{\underline{k}} \partial^{\underline{k}} / {\underline{k}}!$
with $a_{\underline{k}}$ in $\Gamma (\fU, \cO_{\fX })$ satisfying
$\| a_{\underline{k}} \| < c \eta^{\underline{k}}$ for some real constants
$c, \eta$ independent of $\underline{k}$ and $\eta < 1$.

Just as a  reminder let us recall the usual differential
operators. For a smooth formal scheme $\fX$ over
 a complete discrete valuation ring $\V$
let $\I$ be the ideal of the diagonal immersion
$\fX \to \fX \times_{\V} \fX$, the sheaf of principal parts of order $n$ is
the sheaf $\cP^n_{\fX / \V}= \cO_{\fX \times \fX} / \I^{n+1} .$
Considering $\cP^n_{\fX / \V}$ as a $\cO_{\fX}$ module via
$f \mapsto f \otimes 1$ (for $f$ in $\cO_{\fX}$) the sheaf of differential operators
of order $n$ is $\D_{\fX / \V, n}
:= {\cH}om_{\cO_{\fX}} (\cP^n_{\fX / \V},\cO_{\fX}).$ And
the sheaf of differential operators
is $\D_{\fX / \V} := \cup_{n \geq 0} \D_{\fX / \V, n}.$
(This given in EGA IV, 16.8.)

\section{Siegel modular varieties}

\subsection{}

I would like to begin with a sketch of the moduli problem.

Fix integers $g \geq 2$, $ N \geq 3$, a primitive $N$th root of
unity $\zeta_N$.
 Let $Sch$ be the category
of locally noetherian schemes over $\ZZ [1/N, \zeta_N]$.

Consider the functor which assigns to a scheme $S$ in $Sch$ the set
of isomorphism classes of triples $(A/S, \lambda, \alpha)$ where

1) $A \to S$ is an abelian scheme of relative dimension $g$.

2) $\lambda$ is a principal polarization  $A/S \to \hat{A}/S$,
 where $\hat{A}/S$ denotes the dual abelian scheme of
$A/S$.

3) $\alpha$ is a symplectic isomorphism $(A[N], Weil) \to ((\ZZ / N
\ZZ)^g \times \mu_N^g , standard)$ (we call this a level $N$
structure).

\noindent This functor is representable by a quasiprojective scheme
$Y$ smooth over $\ZZ [1/N]$ (Mumford \cite{Mum 94}, \cite{Cha 85}
1.7).

We choose a smooth toroidal compactification ${X}$ of $Y$ such that
$X$ is a smooth projective scheme over Spec $\ZZ [1/N, \zeta_N]$
together with an action of $Sp_{2g} (\ZZ [1/N, \zeta_N)$
(\cite{Ash 75}, \cite{FaC 90} V section5).

Fix a prime $p$ not dividing $N$. Let $\V$ be the  completion of
$\ZZ [1/N, \zeta_N]$ at a prime $\fp$ above $p$. Let $\FF$ (resp. $F$) be
the residue field (resp. field of fractions) of $\V$. Let ${\fX}$ be
the formal completion of $X$ along the closed fibre over $\fp$.

\subsection{}
Let $G$ be the semisimple split symplectic group scheme of rank $g$
over $\ZZ$. We can take
$$G(\ZZ) = \{ \gamma \in GL_{2g}(\ZZ) : {}^{t}\gamma J \gamma =  \; \}
$$
with $J = \left( \begin{array}{ll}
0 & I_{g} \\
-I_{g} & 0
\end{array}\right)$. Write $\Gamma$ for the congruence subgroup of level $N$,
i.e. $\Gamma = \ke (G(\ZZ) \to G(\ZZ / N \ZZ))$. Let $K$ be the
maximal compact subgroup of $G(\RR)$. Then the complex points of $Y$
is
\[ Y(\CC) = K \backslash G(\RR) / \Gamma \; .   \]
In fact let $\fG$ be the formal completion of $G$ at $\fp$. Then we
have an action of $\fG$ on $\fX$. Suppose we have a system of
differential equations or better an
 $\M$  in the bounded derived category
$D_{\fG}^b ({\D}_{\fX}^{\+})$  of $\fG$-equivariant
${\D}_{\fX}^{\+}$ modules, then on the solution spaces $Ext^*(\M,
\cO_{\fX})$ the $p$-adic symplectic group $G(F)$ acts.

\subsection{}

Write
$X_{\FF}$ for the closed fibre of $\fX$. Let $U$ be the ordinary
locus in $Y_{\FF}$. We have a frame $ U \xrightarrow{j} X_{\FF} \to
\fX$.
Let $\M$ be in the bounded derived category $D_{\fG}^b
({\D}_{\fX}^{\+})$  of $\fG$-equivariant ${\D}_{\fX}^{\+}$ modules.
The problem is to study the admissibility of the
$Sp_{2g}(\QQ_p)$ representations on the
$Ext^*(\M, \cO_{\X})$ obtained from the overconvergent solutions
$R{\cH}om_{{\D}_{\fX}^{\+}} (\M,
sp_{*} j^{\+} \cO_{\X} )_F$.

As the Galois group $Gal(\bar{F} /F)$ also acts on $Ext^*(\M, \cO_{\X})$ via
its action on $\X$ as a moduli space, it would be interesting to study
$Ext^*(\M, \cO_{\X})$ as a $Gal(\bar{F} /F) \times Sp_{2g}(\QQ_p)$ module;
in particular the relation with $p$ adic Hodge theory via the Beerger-Kisin functor
between the $(\phi, N)$-modules over $K_0$ and the
$(\phi, N_{\triangledown})$-modules over $\cO$ as in the original work of Berger-Breuil.

\subsection{}
Write $\pi : A \to Y$ for the universal abelian scheme. Put $\omega$
= $\omega_{Y/A}$ = $\pi_* (\wedge^g \Omega^1_{A/Y})$. The
 elements of $H^0 (Y, \omega_{A/Y}^{\otimes k})$
are Siegel modular forms of weight $k$ and level $N$.

The universal abelian scheme $A \to Y$ extends to a semiabelian
scheme  $\bar{A}  \to X$. Write $\bar{\omega}$ for the extension of
$\omega$ to $X$ (\cite{Mum 77}). According to Koecher's principle we
have
\begin{displaymath}
H^0(Y, \omega^{\otimes k}) = H^0(X, \bar{\omega}^{\otimes k}).
\end{displaymath}

Let ${\X}$ be the rigid analytic space associated to ${\fX}$. 
We shall use the same notation for the sheaf on ${\X}$ associated to
the sheaf $\bar{\omega}$ on $X$. We can twist the  $(\GG_m /
K)$-torsor $\bar{\omega}$ on ${\X}$ on ${\X}$ by pushing forward
along a a character $\kappa : \GG_m / K \to  \GG_m / K$ and we
denote the resulting sheaf by $\bar{\omega}^{\kappa}$.

We now apply Berthelot's construction to obtain $j^{\dagger}
\bar{\omega}^{\kappa}$. The rigid analytic space on $Y \otimes K$ is
denoted by $\Y$.  We call an element of $M_{\kappa}(N)^{\dagger} =$
$H^0(\Y, j^{\+} \bar{\omega}^{\kappa})$ an overconvergent Siegel
modular form of weight $\kappa$ and level $N$.

It would be nice to know that $j^{\dagger}\bar{\omega}^{\kappa}$ is
independent of the choice of the lattice in the polyhedral
decomposition used for the construction of the integral model of the toroidal
compactification. As it is we shall simply take the direct limit
of $j^{\dagger}\bar{\omega}^{\kappa}$ over compactifications build from a tower of
lattices to allow the actions of Hecke operators as constructed in
\cite{FaC 90} (see also the theis of KeiWen Lan).

We have a natural $\D_{\fX}^{\+}$-module structure on $sp_{*}
j^{\dagger} \bar{\omega}^{\kappa}$. Are the representations of the
$p$-adic symplectic group on the solution spaces $Ext^*(sp_{*}
j^{\dagger} \bar{\omega}^{\kappa}, \cO_{\X})$ admissible?

We can also ask if we take as $M$ the sheaf of
overconvergent $p$-adic Siegel modular forms or  $j^{\+}$ of
automorphic bundles (\cite{Har 856}, \cite{Mil 88} )
 what do we get?

\subsection{}
I agree that it would be easier to begin with the rank one case
namely that of $SL_2$. It would be fantastic to establish Breuil's
conjecture on $p$ adic Langlands correspondence for $GL_2$ via
arithmetic $\D$-modules.

On the other hand it has also been suggested that one could look at the
universal deformation space of one dimensional formal group of height $h$ together
with its period map into the $h-1$ dimensional projective space.
But such formal schemes are neither of finite type nor proper over the base scheme
and so the standard assumptions in Berthelot are not satisfied. May be there is some way to overcome this.

\medskip
\subsection{}
Though it is said in a different setting I would still like to
recall Varadarajan's comment in the real case: "...resort to the
hathayoga of special functions to do everything...will be an
exercise in futility for it will tell us almost nothing of what is
likely to happen in the general case." (\cite{VSV 89} p.226.) In
another context in which Katz studied $p$-adic $L$ functions by
applying $p$-adic differential operators to Eisenstein series, he
wrote : "we were able to prove theorems by dipping into classical
material... We were for a long time blinded by these riches to the
simple cohomological mechanism which in some sense underlies them."
(\cite{Kat 78} p.205.)

\bigskip
\bigskip

\end{document}